\documentclass[a4paper,11pt]{article}
\usepackage[centertags]{amsmath}
\usepackage{amssymb}
\usepackage{amsthm}
\usepackage{graphicx}
\usepackage[pagewise]{lineno}
\usepackage{amsfonts}
\usepackage{mathrsfs}
\usepackage{times,cite}
\usepackage{geometry}
\usepackage{color}
\newtheorem{thm}{Theorem}[section]
\newtheorem{lem}[thm]{Lemma}
\newtheorem{prop}[thm]{Proposition}
\newtheorem{cor}[thm]{Corollary}
\newtheorem{defn}[thm]{Definition}
\theoremstyle{definition}
\newtheorem{rem}[thm]{Remark}

\newcommand{\e}{\varepsilon}
\newcommand{\p}{\partial}

\newenvironment {Proof} {\noindent {\bf Proof. }}{\hfill $\square$\par\vspace{3mm}}
\newenvironment {ProofTheorem} {\noindent {\bf Proof of Theorem \ref{thm2}}}{\hfill $\square$\par\vspace{3mm}}
\numberwithin{equation}{section}
\allowdisplaybreaks    

\newenvironment{acknowledgment}{\noindent{\bf Acknowledgment }}{}
\date{}
\title{Global well-posedness of the Navier-Stokes equations with Navier-slip boundary conditions in a strip domain}
\author{{ Quanrong Li$^a$, Shijin Ding$^b$ \thanks{Corresponding author. Emails: {\it quanrong\_li@szu.edu.cn(Q. Li),dingsj@scnu.edu.cn(S. Ding)}}}\\
{\it\small $^a$College of Mathematics and Statistics,} {\it\small Shenzhen University,}\\
{\it\small Shenzhen, 518060, Guangdong, China}\\
{\it\small $^b$South China Research Center for Applied Mathematics and Interdisciplinary Studies,}\\
{\it\small South China Normal University,}{\it\small Guangzhou, 510631, Guangdong, China}\\}
\begin{document}
\newcommand{\D}{\displaystyle}
\maketitle

 \begin{abstract}
 This paper is concerned with the existence and uniqueness of the strong solution to the incompressible Navier-Stokes equations with Navier-slip boundary conditions in a two-dimensional strip domain where the slip coefficients may not have defined sign. In the meantime, we also establish a number of Gagliardo-Nirenberg inequalities in the corresponding Sobolev spaces which will be applicable to other similar situations.
 \end{abstract}

{\bf Keywords: }Navier-Stokes equations, Global well-posedness, Navier-slip conditions.

{\em AMS Subject Classification:} 76N10, 35Q30, 35Q35.
\section{Introduction}

 Navier-Stokes equations is one of the most classical mathematical models in fluid dynamics and is also the basic system in the study of most complex fluids. Since being derived by the famous physicists C. Navier and G. Stokes, it has attracted the attentions of considerable number of mathematicians and physicists during the past over 100 years. Precisely, the incompressible Navier-Stokes equation reads as follows
\begin{equation}\label{a.1}
\begin{cases}
\p_t\rho+\mathrm{div}(\rho\mathbf{v})=0,\\
\p_t(\rho\mathbf{v})+\mathrm{div}(\rho\mathbf{v}\otimes\mathbf{v})+\nabla p=\mathrm{div}\mathbb{S}+\rho\mathbf{f},\\
\p_t(\rho E)+\mathrm{div}(\rho\mathbf{v}E+\mathbf{v}p)=\mathrm{div}(\mathbf{v}\mathbb{S}+\kappa\nabla\theta)+\rho\mathbf{f}\cdot\mathbf{v},\\
\mathrm{div}\mathbf{v}=0,
\end{cases}
~~\textrm{in~~}~\Omega\times(0,T).
\end{equation}
where $\rho,\mathbf{v}, \theta, E$ are density, velocity, absolute temperature and total energy of the fluid, respectively and $p, \mathbb{S}, \mathbf{f}$ stand for the pressure, stress tensor and external force, respectively. We point out that the first three equations are deduced by the conservation of {\it mass}, {\it momentum} and {\it energy}, respectively. For the derivation of the Navier-Stokes system, we refer the reader to books by G. P. Galdi\cite{Galdi}, by R. Temam\cite{Temam84} and by P. L. Lions\cite{PLions1}.

If the density and temperature are constants and the tress tensor is taken to be the simplest one $\mathbb{S}=\mu(\nabla\mathbf{v}+\nabla^T\mathbf{v})$, then system (\ref{a.1}) reduces to the following
\begin{equation}\label{a.2}
\begin{cases}
\p_t\mathbf{v}+\mathbf{v}\cdot\nabla\mathbf{v}+\nabla p=\mu\Delta\mathbf{v}+\mathbf{f},\\
\mathrm{div}\mathbf{v}=0,
\end{cases}
\ \textrm{in~~}\Omega\times(0,T).
\end{equation}
To seek solutions for \eqref{a.2} and study the properties of the solutions, it is necessary to impose some conditions, such as $\Omega$ is a bounded domain in $\mathbb{R}^d,d\geq 2$, and the velocity satisfies the Dirichlet boundary condition $\mathbf{v}|_{\p\Omega}=\varphi$ and initial conditions $\mathbf{v}(0)=\mathbf{v}_0$; or $\Omega=\mathbb{R}^d,d\geq 2$, then give the data of $\mathbf{v}$ in the far field and the initial time, which is called Cauchy problem. In all these cases, to our knowledge, the uniqueness of the weak solution to the system $(\ref{a.2})$ in 3D with general initial data $\mathbf{v}_0$, or equivalently, the higher order regularity of the weak solution, is still an open problem.

It should be noted that most of the existing results mainly focus on the Dirichlet boundary conditions, i.e.no-slip boundary conditions. However, there are many other kinds of boundary conditions which also match with the reality. For example, hurricanes and tornadoes do slip along the ground and lose energy as they slip\cite{Olive07}. In 1827, the famous mathematician and physicist C. Navier\cite{CNavier} first considered the slip phenomena and proposed the following boundary conditions, called Navier-slip boundary conditions:
\begin{equation}\label{a.3}
\begin{cases}
\mathbf{v}\cdot\mathbf{n}=0,\\
2\mu\mathbb{D}(\mathbf{v})\mathbf{n}\cdot\tau=k(\mathrm{x})\mathbf{v}\cdot\tau,
\end{cases}
\ \textrm{on~~}\ \p\Omega,
\end{equation}
where $\mathbb{D}(\mathbf{v})=\frac{1}{2}(\nabla\mathbf{v}+\nabla^T\mathbf{v})$, $ \mathbf{n}$ and $\tau$ are unit outer normal vector and tangential vector of the boundary $ \p\Omega$. In (\ref{a.3}), $ k(\mathrm{x})$ is a physical parameter, which can be a constant, function in $ L^\infty(\partial\Omega) $\cite{Kelliher06} and even a smooth metrix\cite{Kelliher12}. Here we consider the case that $k(\mathrm{x})$ is constant, called the slip coefficient.

We should also mention that the most known discussions on the Navier-slip boundary value problems are for the "classical" cases in which the slip coefficients are non-positive, that is, $k(\mathrm{x})\leq 0$ in the boundary conditions \eqref{a.3}, which is according with the friction effect. The pioneers in analysing the Navier-Stokes equations with Navier-slip boundary conditions should be Solonnikov and \v{S}\v{c}adilov \cite{Solonnikov}, who considered the linearized equations in steady case. Afterwards, B. da Veiga \cite{daVeiga} established the existence and the regularity of the weak solutions for the nonlinear problem in the upper half space, while C. Amrouche et al.\cite{Amrouche11,Amrouche14} gave the corresponding results in bounded domain and external domain.

What we are interested in this paper is for the "non-classical" cases in which the slip coefficients may be positive, and the domain is unbounded. As being pointed out by Serrin\cite{Serrin} in 1959, $k(\mathrm{x})$ is unnecessary to be negative. Moreover, there do be some phenomena in the real world with $k(\mathrm{x})>0$. For example, the effective slip length $\alpha$ on the flat gas-liquid interface is always positive\cite{Haase}. Navier-slip boundary conditions (\ref{a.3}) with $k(\mathrm{x})>0$ is also widely applied in the numerical modeling of fluid with rough boundary, such as in aeronautical dynamics or in the permeable boundary, where (\ref{a.3}) are called Beavers-Joseph law\cite{Joseph67,Amrouche11}, in the weather forecast and Hemodynamics\cite{Joseph67,Chauhan93}), or some case where the boundary accelerates the fluid\cite{Accelerate,Olive07}. The readers could refer to Y. Xiao et al.\cite{YXiao07,YXiao13} and the reference therein for some results on the vanishing viscosity limit of the time-depending Navier-Stokes equtions. For more details in physical applications and numerical analysis, please refer to\cite{Achdou,Bansch01,Joseph67,Jager01,Jager00,VJohn,TQian,Serrin}.

In 2017, H. Li and X. Zhang\cite{lizhang} established the global well-posednes of the 3D compressible Navier-Stokes equations in a strip domain with Dirichlet boundary condition on the upper plane and Navier-slip boundary condition on the bottom. However, the slip coefficient must be an negative constant. In 2018,  Xin and the authors of this paper published a paper \cite{DL18} on the stability analysis for the Navier-slip boundary value problems for this "non-classical" cases. We found in \cite{DL18} that if some of the slip coefficients are positive, the kinetic energy generated on the boundary may cause instability if the viscosity is not large enough. So, we defined, in \cite{DL18},  a critical viscosity expressed only by the slip coefficients to distinguish the stability from the instability. However, in that paper we did not give the detailed proof for the global existence and uniqueness of strong solutions for the so called "non-classical" cases. We find that the proof of the global well-posedness is not a trivial problem mainly because of the unboundedness of the domain and the boundary conditions. Moreover, in this paper, we have also derived a number of Gagliardo-Nirenberg inequalities which may be applicable to other similar cases.

The rest of this paper will be arranged as follows: in Section 2 , some notations will be given and the definition together with the main theorems will be stated; Section 3 is arranged for preliminary, that is, the proof of a series of Gagliardo-Nirenberg inequality; the global existence of the unique weak solution to system \eqref{b.1} will be established in Section 4 and the proof of higher order regularity to the weak solution, so that the weak solution is in fact a strong one, will be given in Section 5.

\section{Notations and main results}

Precisely, we consider the following initial boundary value problem
\begin{equation}\label{b.1}
\begin{cases}
\p_t\mathbf{u}+\mathbf{u}\cdot\nabla\mathbf{u}+\nabla p=\mu\Delta\mathbf{u}+\mathbf{f},&~\textrm{~in~~}\Omega\times(0,T),\\
\nabla\cdot\mathbf{u}=0,&~\textrm{~in~~}\Omega\times(0,T),\\
\mathbf{u}\cdot\mathbf{n}=0,&~\textrm{~on~~}\p\Omega\times(0,T)\\
2\mu\mathbb{D}(\mathbf{u})\mathbf{n}\cdot\tau=k(x,y)\mathbf{u}\cdot\tau,&~\textrm{~on~~}\p\Omega\times(0,T)\\
\mathbf{u}(0)=\mathbf{u}_0,&~\textrm{~in~~}\Omega.
\end{cases}
\end{equation}
where $ k(x,1)=k_1 $,$k(x,0)=k_0$ are constants and $\Omega:=\mathbb{R}\times(0,1)$. For convenience, we denote
\begin{align*}
&\mathbf{H}:=\{\mathbf{v}\in L^2|\nabla\cdot\mathbf{v}=0,\ \mathbf{v}\cdot\mathbf{n}=0, \textrm{~on~}\p\Omega\},\\
&\mathbf{V}:=\{\mathbf{v}\in H^1|\nabla\cdot\mathbf{v}=0,\ \mathbf{v}\cdot\mathbf{n}=0, \textrm{~on~}\p\Omega\},\\
&\mathcal{W}:=\{\mathbf{v}\in \mathbf{V}\cap H^2|\mathbf{v}\textrm{~~satisfies~~(\ref{a.3})}_2\}.
\end{align*}
 In the meantime, we denote $L^2(0,1)$ and $H^k(0,1)$ by $L^2$ and $H^k$, for simplicity. Without confusion, we will also write $L^p(\Omega)$ and $H^k(\Omega)$  by $L^p$ and $H^k$, respectively. The integral form $ \int_\Omega fdxdy$ will be simply denoted by $\int f$. In addition, the scalar function and vector function will be denoted by $f$ and $\mathbf{f}$ for distinction, such as $\mathbf{f}=(f^1,f^2) $, but the product functional space $(X)^2$ will also be denoted by $X $. For example, the vector function $\mathbf{u}\in (H^1)^2$ will be still denoted by $\mathbf{u}\in H^1$. The usual notations will be used as in general unless extra statement.

We will first prove the global existence of the unique weak solution to (\ref{b.1}), and then improve the regularity to reach the global strong solution. Now we give the definition of weak solutions.

\begin{defn}\label{def}
$\mathbf{u}$ is a weak solution to the initial boundary value problem (\ref{b.1}) defined in $\Omega\times(0,T)$, if it satisfies

1. $\mathbf{u}\in L^\infty(0,T;\mathbf{H})\cap L^2(0,T;\mathbf{V})$; 2. $\mathbf{u}(x,y,0)=\mathbf{u}_0(x,y)$, {\mbox{in}} \ $ \Omega$;

3. For any $\mathbf{v}\in\mathbf{V}$, there holds that
\begin{align*}
&\frac{d}{dt}\int_{\Omega}\mathbf{u}(t)\cdot\mathbf{v}+\mu\int_{\Omega}\mathbb{D}(\mathbf{u}(t)):\mathbb{D}(\mathbf{v})
+\int_{\Omega}\mathbf{u}(t)\cdot\nabla\mathbf{u}(t)\cdot\mathbf{v}\\
&=\int_{\p\Omega}k(x,y)\mathbf{u}(t)\cdot\mathbf{v}dS+\int_{\Omega}\mathbf{f}\cdot\mathbf{v}.
\end{align*}
\end{defn}

For the first step, we prove the following global well-posedness in weak sense
\begin{thm}\label{thm1}
For any initial data $\mathbf{u}_0\in\mathcal{W}$ and external force $\mathbf{f}\in H^1(0,T;L^2)$, there exists an unique weak solution $\mathbf{u}\in L^\infty(0,T;\mathbf{V})\cap L^2(0,T;\mathbf{V})$ to the initial boundary value problem (\ref{b.1}), which satisfies $\p_t\mathbf{u}\in L^\infty(0,T;\mathbf{H})\cap L^2(0,T;\mathbf{V})$ and
\begin{align}\label{c.2}
&\sup_{0\leq t\leq T}\left(\|\mathbf{u}(t)\|^2_{H^1}+\|\p_t\mathbf{u}(t)\|^2_{L^2}\right)
+\int_0^T\left(\|\mathbf{u}(t)\|^2_{H^1}+\|\p_t\mathbf{u}(t)\|^2_{H^1}\right)\nonumber\\
&\leq C\left(T,\|\mathbf{u}_0\|_{H^2},\|\mathbf{f}\|_{H^1(0,T;L^2)}\right).
\end{align}
\end{thm}

Then, the main result of this paper is to prove the following global well-posedess theory

\begin{thm}\label{thm2}
For any $\mathbf{u}_0\in\mathcal{W},T>0 $ and $\mathbf{f}\in H^1(0,T;L^2) $ initial boundary value problem (\ref{b.1})exists an unique strong solution $\mathbf{u}\in L^\infty(0,T;H^2)$ satisfies
\begin{align}\label{b.2}
\sup_{0\leq t\leq T}\left(\|\mathbf{u}(t)\|_{H^2}+\|\nabla p(t)\|_{L^2}\right)\leq C\left(T,\|\mathbf{u}_0\|_{H^2},\|\mathbf{f}(t)\|_{H^1(0,T;L^2)}\right).
\end{align}
\end{thm}

Before continuing, we would like to have some words on the main result of this paper.

\begin{rem}\label{rem1}
(1)~~Theorem \ref{thm2} is valid without any smallness constraint on the initial data, which is consistent with the existed results on dimension 2.
(2)~~The global existence of the unique strong solution is valid without any constraint on the sign of the slip coefficients.
\end{rem}
\section{Preliminary}

Note that the domain $\Omega$ is unbounded and the boundary $\p\Omega$ is non-compact, which lead to the difficulty in finding the smooth orthonormal basis for the construction of Galerkin approximate solutions. Thus, we first find solutions in a subdomain $\Omega_L:=(-L,L)\times(0,1)$ with the similar Navier-slip boundary conditions on $\{-L,L\}\times[0,1]$, where $ k(-L,y)=k(L,y)=0$. We infer that the definition of weak solution is similar to that on $\Omega$  and denote the constraint of $ \mathbf{H},\mathbf{V},\mathcal{W}$ in $\Omega_L$ by $\mathbf{H}_L,\mathbf{V}_L,\mathcal{W}_L$, respectively. Without lose of any generality, we take $L\geq 1.$

To apply the Galerkin method in proving the existence of the unique solution, we need the following two lemmas, which are similar to Lemma 2.1 and Lemma 2.2 of \cite{Clopeau}, respectively.

\begin{lem}\label{mc.2}
Assume that $\mathbf{v}\in H^2(\Omega_L)$ satisfies $\mathbf{v}\cdot\mathbf{n}=0$ on the boundary $\p\Omega_L$. Then it holds on the boundary that
\[2\mathbb{D}(\mathbf{v})\mathbf{n}\cdot\tau=\mathrm{curl}\mathbf{v},\]
where $\mathrm{curl}\mathbf{v}:=\p_xv^2-\p_yv^1$.
\end{lem}
\Proof
Refer to the proof of Lemma 2.1 in \cite{Clopeau}.
\endProof
\begin{lem}\label{mc.3}
There exists a basis $\{\mathbf{w}_1,\mathbf{w}_2,\cdots,\mathbf{w}_n,\cdots\}\subset H^3(\Omega_L)$ to $\mathbf{V}_L$, such that
\[2\mu\mathbb{D}(\mathbf{w}_m)\mathbf{n}\cdot\tau=k(x,y)\mathbf{w}_m\cdot\tau,\ \ \textrm{on~~}\p\Omega_L,\ m=1,2,\cdots.\]
Moreover, $\{\mathbf{w}_1,\mathbf{w}_2,\cdots,\mathbf{w}_n,\cdots\}$ is also an orthonormal basis of $\mathbf{H}_L$.
\end{lem}
\Proof
The main idea of the proof, which consists of three steps, is quite different from that of Lemma 2.2 in \cite{Clopeau}, for the domain here is a rectangular region.

{\it Step 1.} For some positive constant $\beta$ large enough, consider the auxiliary eigenvalue problem
\begin{align}\label{c.02}
\begin{cases}
-\mu\Delta\mathbf{u}+\nabla p+\beta\mathbf{u}=\lambda\mathbf{u},&\ \textrm{in~~}\ \Omega_L,\\
\mathrm{div}\mathbf{u}=0,&\ \textrm{in~~}\ \Omega_L,\\
\mathbf{u}\cdot\mathbf{n}=0,&\ \textrm{on~~}\ \p\Omega_L,\\
\mu\mathbb{D}(\mathbf{u})\mathbf{n}\cdot\tau=k(x,y)\mathbf{u}\cdot\tau,&\ \textrm{on~~}\ \p\Omega_L,\\
\end{cases}
\end{align}
of which the variational form is to seek $\mathbf{u}\in V_L$ and  $\lambda\neq 0$, such that for any $\mathbf{v}\in V_L$, there holds
\begin{align}\label{c.03}
\mu\int_{\Omega_L}\nabla\mathbf{u}\cdot\nabla\mathbf{v}+\beta\int_{\Omega_L}\mathbf{u}\cdot\mathbf{v}
+\int_{\p\Omega_L}k(x,y)(\mathbf{u}\cdot\tau)(\mathbf{v}\cdot\tau)=\lambda\int_{\Omega}\mathbf{u}\cdot\mathbf{v}.
\end{align}

Note that the bilinear form
\[a(\mathbf{u},\mathbf{v}):=\mu\int_{\Omega_L}\nabla\mathbf{u}\cdot\nabla\mathbf{v}+\beta\int_{\Omega_L}\mathbf{u}\cdot\mathbf{v}
+\int_{\p\Omega_L}k(x,y)(\mathbf{u}\cdot\tau)(\mathbf{v}\cdot\tau)\]
is continuous and symmetric on $V_L\times V_L$. In particular, when $\mathbf{v}=\mathbf{u}$, one gets for any $\e>0$ that
\begin{align*}
\int_{\p\Omega_L}k(x,y)(\mathbf{u}\cdot\tau)^2&=\int_{\Omega_L}[((k_0+k_1)y-k_0)(u^1)^2]_y\\
&\geq -\e\|\p_yu^1\|^2_{L^2(\Omega_L)}+(k_0+k_1-\e^{-1}\max\{k_1^2,k_0^2\})\|u^1\|^2_{L^2(\Omega_L)},
\end{align*}
Then, as long as $\beta$ being so large that
\[\beta>\beta_0:=\e^{-1}\max\{k_1^2,k_0^2\}-(k_0+k_1),\ \e\in(0,\mu),\]
the bilinear form $a(\cdot,\cdot)$ satisfies
 \[a(\mathbf{u},\mathbf{u})\geq (\mu-\e)\|\nabla\mathbf{u}\|^2_{L^2(\Omega_L)}+(\beta-\beta_0)\|\mathbf{u}\|^2_{L^2(\Omega_L)}
 \geq c_0\|\mathbf{u}\|^2_{H^1(\Omega_L)}.\]
 This indicates that $a(\cdot,\cdot)$ is coercive on $(V_L, V_L)$.

{\it Step 2.} It is clear that the embedding map $V_L\hookrightarrow H_L$ is compact. Then, it follows from the spectral theory of operators that there exists countable eigenvalues $\{\lambda_j\}$ to problem \eqref{c.03} such that as $j\rightarrow +\infty$, $\lambda_j\rightarrow +\infty$. Correspondingly, the eigenfunctions $\{\mathbf{w}_j\}$ constitute a basis of $V_L$, which, in the meantime, is also a orthonormal basis of $H_L$. This means the eigenvalue problem of the Stokes operator
\begin{align}\label{c.04}
\begin{cases}
-\mu\Delta\mathbf{u}+\nabla p=\Lambda\mathbf{u},&\ \textrm{in~~}\ \Omega_L,\\
\mathrm{div}\mathbf{u}=0,&\ \textrm{in~~}\ \Omega_L,\\
\mathbf{u}\cdot\mathbf{n}=0,&\ \textrm{on~~}\ \p\Omega_L,\\
\mu\mathbb{D}(\mathbf{u})\mathbf{n}\cdot\tau=k(x,y)\mathbf{u}\cdot\tau,&\ \textrm{on~~}\ \p\Omega_L,\\
\end{cases}
\end{align}
possesses countable eigenvalues $\{\Lambda_j\}_{j=1}^\infty$, satisfying $-\beta_0<\Lambda_1<\Lambda_2<\cdots$ and $\Lambda_j\rightarrow+\infty$, as $j\rightarrow +\infty$(Ref. L. Evans\cite{Evans}$\S6.2$).

{\it Step 3.} Now, we apply bootstrap method to improve regularity of the eigenfunctions $\{\mathbf{w}_j\}_{j=1}^\infty$. As $\mathbf{w}_j$ satisfies divergence free condition $\mathrm{div}\mathbf{w}_j=0$, there exists steam function $\psi_j$ such that $\mathbf{w}_j=(-\p_y,\p_x)\psi_j$. Further, denote $\omega_j:=\p_x\mathrm{w}_j^2-\p_y\mathrm{w}_j^1$. Then $\psi_j$ satisfies the following Dirichlet problem
\begin{align}\label{c.05}
\begin{cases}
-\Delta\psi_j=-\omega_j,&\ \textrm{in~~}\ \Omega_L,\\
\psi_j=0,&\ \textrm{on~~}\ \p\Omega_L.
\end{cases}
\end{align}
In virtue of \eqref{c.04} together with Lemma \ref{mc.2}, we deduce that $W_j=\omega_j-g$ satisfies the Dirichlet boundary value problem
\begin{align}\label{c.06}
\begin{cases}
-\mu\Delta W_j=\lambda_j \omega_j+\mu\Delta g,&\ \textrm{in~~}\ \Omega_L,\\
W_j=0,&\ \textrm{on~~}\ \p\Omega_L,
\end{cases}
\end{align}
where $g(x,y):=[(k_0+k_1)y-k_0]\mathrm{w}^1(x,y)$.

 Note that $\mathbf{w}_j\in H^1(\Omega_L)$, namely, the right-hand side of \eqref{c.06}$_1$ belongs to $H^{-1}(\Omega_L)$. Then it follows from the elliptic estimate that $W_j\in H^1_0(\Omega_L)$, which further implies that $\omega_j\in H^1(\Omega_L)$. Consequently, applying the theory of elliptic equations to system \eqref{c.05} yields $\psi_j\in H^3(\Omega_L)$, which indicates $\mathbf{w}_j\in H^2(\Omega_L)$. The proof of this lemma is completed.
\endProof

In general, the uniform constants in Gagliardo-Nirenberg inequalities depend on the shape or the size of the domain. To obtain the energy estimates independent of $L$, we need the following Gagliardo-Nirenberg inequalities, of which the uniform constants are independent of the horizontal length of the rectangular domain. The authors believe that these inequalities will be applicable in other similar situations.
\begin{lem}\label{mc.4}
{\bf( $L^2(\Omega_L)$ estimate)}~~There exists a constant $C>0$, being independent of $L$, such that for any $\mathbf{u}\in\mathbf{V}_L$, there holds
\begin{align}\label{c.3}
\|\mathbf{u}\|_{L^2(\Omega_L)}\leq C\|\p_y\mathbf{u}\|_{L^2(\Omega_L)}.
\end{align}
\end{lem}

\Proof First prove for $u^2$. Since $u^2(x,0)=u^2(x,1)=0$, there holds $u^2(x,y)=\int_0^yu^2_y(x,\theta)d\theta$. Then
\[\int_{\Omega_L}|u^2(x,y)|^2dxdy=\int_{\Omega_L}\left|\int_0^yu^2_y(x,\theta)d\theta\right|^2dxdy\leq \int_{\Omega_L}|u_y^2(x,y)|^2dxdy.\]

To $u^1$, it follows from the incompressible condition that $\int_0^1u^1_x(x,y)dy=-\int_0^1u^2_y(x,y)dy=0$, i.e. $\int_0^1u^1(x,y)dy$ is a constant. Besides, it is clear that
\[0=\int_{\p\Omega_L}x\mathbf{u}\cdot\mathbf{n}dS=\int_{\Omega_L}\mathrm{div}(x\mathbf{u})=\int_{\Omega_L}u^1(x,y)dxdy,\]
which means $\int_0^1u^1(x,y)dy\equiv 0$. The it follows from the Poincar\'e inequality on the vertical direction that
 \[\int_0^1|u^1(x,y)|^2dy\leq C\int_0^1|u^1_y(x,y)|^2dy\]
holds for some constant $C>0$. Integrating this inequality respect to $x$ completes the proof of this lemma.
\endProof

\begin{cor}\label{cc.5}
{\bf($L^2(\Omega)$ estimate)}~~There exists a constant $C>0$, such that for any $\mathbf{u}\in\mathbf{V}$, it is valid that
\begin{align}\label{c.4}
\|\mathbf{u}\|_{L^2}\leq C\|\p_y\mathbf{u}\|_{L^2}.
\end{align}
\end{cor}

\begin{lem}\label{mc.6}
{\bf ($L^4(\Omega_L)$ estimate)}~~There exists a constant $C>0$, being independent of $L$, such that for any $\mathbf{u}\in\mathbf{W}_L$, there holds
\begin{align}\label{c.5}
\|\mathbf{u}\|^2_{L^4(\Omega_L)}\leq C\|\mathbf{u}\|_{L^2(\Omega_L)}\|\nabla\mathbf{u}\|_{L^2(\Omega_L)}.
\end{align}
\end{lem}

\Proof Note that the boundary conditions (\ref{b.1})$_3$ and (\ref{b.1})$_4$ can be rewritten as
\begin{align*}
&u^1(\pm L,y)=0;\ u^2(x,0)=u^2(x,1)=0;\\
&\mu\p_yu^1(x,0)=-k_0u^1(x,0),\mu\p_yu^1(x,1)=k_1u^1(x,1);\p_xu^2(\pm L,y)=0.
\end{align*}
We first claim that if $f\in H^1(\Omega_L)$ satisfies $f(-L, y)=f(x,0)=0$, then
\begin{align}\label{c.01}
\|f\|^2_{L^4}\leq 2\|f\|_{L^2}\|\nabla f\|_{L^2}.
\end{align}
In fact, we have
\[|f(x,y)|^2=2\int_{-L}^x f(s,y)f_x(s,y)ds\leq 2\|f(y)\|_{L^2(-L,L)}\|f_x(y)\|_{L^2(-L,L)},\]
 and
\[|f(x,y)|^2=2\int_0^y f(x,\theta)f_y(x,\theta)d\theta\leq 2\|f(x)\|_{L^2(0,1)}\|f_y(x)\|_{L^2(0,1)}.\]
Multiplying the above two equations, integrating over $\Omega_L$, and using H\"older inequality, we get
\begin{align*}
\int_{\Omega_L}|f(x,y)|^4&\leq 4\int_0^1\|f(y)\|_{L^2(-L,L)}\|f_x(y)\|_{L^2(-L,L)}
\int_{-L}^L\|f(x)\|_{L^2(0,1)}\|f_y(x)\|_{L^2(0,1)}\\
&\leq 4\|f\|^2_{L^2(\Omega_L)}\|f_x\|_{L^2(\Omega_L)}\|f_y\|_{L^2(\Omega_L)}.
\end{align*}
Then (\ref{c.01}) follows.

Now, we prove (\ref{c.5}) for $u^1$. Denote $\zeta(y)$ to be a smooth cut-off function on $\mathbb{R}$ satisfies (1) when $|y|\leq 1$, $\zeta(y)\equiv 1$; (2) when $|y|\geq 2$, $\zeta(y)\equiv 0$; (3) for any $y\in\mathbb{R}$ there holds $ \zeta(y)\in[0,1]$ and $|\zeta^\prime(y)|\leq 2$; (4) $\zeta(y)=\zeta(-y)$.
Further denote $v^1(x,y):=\zeta(2y)\tilde{u}^1(x,y)$, where
\begin{equation*}
\tilde{u}^1(x,y):=
\begin{cases}
e^{\frac{k_0}{\mu}y}u^1(x,y),&\ y\in [0,1],\\
e^{\frac{-k_0}{\mu}y}u^1(x,-y),&\ y\in [-1,0].
\end{cases}
\end{equation*}
Then $v^1\in H^1([-L,L]\times[-1,1])$ with $v^1(x,-1)=v^1(x,1)=v^1(\pm L,y)=0$. It can be deduced by (\ref{c.01}) together with the symmetry of $v^1$ that
\[\left(\int_{\Omega_L}|v^1(x,y)|^4\right)^{1/2}\leq 2\sqrt{2}\left(\int_{\Omega_L}|v^1(x,y)|^2\right)^{1/2}\left(\int_{\Omega_L}|\nabla v^1(x,y)|^2\right)^{1/2}.\]
In virtue of the definitions of $v^1$ and $\zeta$, one can rewrite the above inequality as
\begin{align*}
\left(\int_{-L}^L\int_0^\frac{1}{2}|u^1(x,y)|^4\right)^{1/2}\leq C(k_0,\mu)\left(\|u^1\|_{L^2(\Omega_L)}
\|\nabla u^1\|_{L^2(\Omega_L)}+\|u^1\|_{L^2(\Omega_L)}^2\right)\end{align*}
Similarly, it is valid that
\begin{align*}
\left(\int_{-L}^L\int_\frac{1}{2}^1|u^1(x,y)|^4\right)^{1/2}\leq C(k_1,\mu)\left(\|u^1\|_{L^2(\Omega_L)}
\|\nabla u^1\|_{L^2(\Omega_L)}+\|u^1\|_{L^2(\Omega_L)}^2\right)\end{align*}
Adding them up and using (\ref{c.3}) yield inequality (\ref{c.5}) for $u^1$.

In what follows, we prove inequality (\ref{c.5}) for $u^2$. Write $v^2:=\zeta(\frac{x}{L}+1)\tilde{u}^2$ with
\begin{equation*}
\tilde{u}^2(x,y):=
\begin{cases}
u^2(x,y),&\ x\in [-L,L]\\
u^2(-x-2L,y),&\ x\in [-3L,-L].
\end{cases}
\end{equation*}
Then $v^2(x,y)$ satisfies $v^2(-3L,y)=v^2(L,y)=v^2(x,0)=v^2(x,1)=0$.

By (\ref{c.01}) and the symmetry of $v^2$, we have
\[\left(\int_{\Omega_L}|v^2(x,y)|^4\right)^{1/2}\leq 2\sqrt{2}\|v^2(x,y)\|_{L^2(\Omega_L)}\|\nabla v^2(x,y)\|_{L^2(\Omega_L)}.\]
Similarly, according to the definition of $v^2,\zeta$, we rewrite the above inequality as
\begin{align*}
\left(\int_{-L}^0\int_0^1|u^2(x,y)|^4\right)^{1/2}\leq C\left(\|u^2\|_{L^2(\Omega_L)}
\|\nabla u^2\|_{L^2(\Omega_L)}+\|u^2\|_{L^2(\Omega_L)}^2\right).
\end{align*}
We should point out here that the constant $C$ depends on $L^{-1}$ because of the derivation of $\zeta(\frac{x}{L}+1)$. However, since our final end is to take $L\rightarrow+\infty$, so, without lose of any generality, we take $L\geq 1$, and then $C$ is independent of $L$.

Besides, we also have
\begin{align*}
\left(\int_0^L\int_0^1|u^1(x,y)|^4\right)^{1/2}\leq C\left(\|u^2\|_{L^2(\Omega_L)}
\|\nabla u^2\|_{L^2(\Omega_L)}+\|u^2\|_{L^2(\Omega_L)}^2\right)
\end{align*}
Adding them up and using (\ref{c.3})deduce inequality (\ref{c.5}) for $u^2$.
\endProof

As the result of Lemma \ref{mc.6} is independent of $L$, we take $L\rightarrow\infty$ and yield the desired $ L^4$ estimate for functions in $\mathcal{W}$ as follows.

\begin{cor}\label{cc.7}
{\bf($L^4(\Omega)$ estimate)}~~There exists a constant $C>0$, such that for any $\mathbf{u}\in\mathcal{W}$, it holds that
\begin{align}\label{c.6}
\|\mathbf{u}\|^2_{L^4}\leq C\|\mathbf{u}\|_{L^2}\|\nabla\mathbf{u}\|_{L^2}.
\end{align}
\end{cor}

\begin{lem}\label{mc.8}
{\bf($L^2(\Omega_L)$ estimate for gradient)}~~There exists constant $C>0$, being independent of $L$, such that for any $\mathbf{u}\in\mathcal{W}_L$, there holds
\begin{align}\label{c.7}
\|\nabla\mathbf{u}\|^2_{L^2(\Omega_L)}\leq C\|\mathbf{u}\|_{L^2(\Omega_L)}\|\nabla\mathbf{u}\|_{H^1(\Omega_L)}.
\end{align}
\end{lem}

\Proof
In fact, in virtue of integrating by parts and H\"older inequality, we have
\begin{align*}
\mu\int_{\Omega_L}|\nabla\mathbf{u}|^2&=-\mu\int_{\Omega_L}\mathbf{u}\cdot\Delta\mathbf{u}+k_1\int_{-L}^L|u^1(x,1)|^2
+k_0\int_{-L}^L|u^1(x,0)|^2\\
&\leq \mu\int_{\Omega_L}|\mathbf{u}||\Delta\mathbf{u}|+\int_{\Omega_L}[((k_1+k_0)y-k_0)|u^1(x,y)|^2]_ydxdy\\
&\leq \mu\|\mathbf{u}\|_{L^2(\Omega_L)}\|\Delta\mathbf{u}\|_{L^2(\Omega_L)}+C\|\mathbf{u}\|^2_{L^2(\Omega_L)}
+C\|\mathbf{u}\|_{L^2(\Omega_L)}\|\nabla\mathbf{u}\|_{L^2(\Omega_L)}\\
&\leq C\|\mathbf{u}\|_{L^2(\Omega_L)}\|\nabla\mathbf{u}\|_{H^1(\Omega_L)},
\end{align*}
in which (\ref{c.3}) has been used in the last inequality.
\endProof

\begin{cor}\label{cc.9}
{\bf($L^2(\Omega)$ estimate for gradient)}~~There exists a constant $C>0$, such that foe any $\mathbf{u}\in\mathcal{W}$, it is true that
\begin{align}\label{c.8}
\|\nabla\mathbf{u}\|^2_{L^2}\leq C\|\mathbf{u}\|_{L^2}\|\nabla\mathbf{u}\|_{H^1}.
\end{align}
\end{cor}

\begin{lem}\label{mc.10}
{\bf(Korn's inequality on $\Omega_L$)}~~There exists a constant $C>0$, being independent of $L$, such that for any $\mathbf{u}\in \mathbf{V}_L$, it holds that
\begin{align}\label{c.9}
\|\mathbb{D}(\mathbf{u})\|_{L^2(\Omega_L)}\geq C\|\mathbf{u}\|_{H^1(\Omega_L)}.
\end{align}
\end{lem}

\Proof
Note that
\begin{align}\label{c.10}
\int_{\Omega_L}|\nabla\mathbf{u}+\nabla^T\mathbf{u}|^2=2\int_{\Omega_L}|\nabla\mathbf{u}|^2
+2\sum_{i,j=1}^2\int_{\Omega_L}\p_iu^j\p_ju^i.
\end{align}
In virtue of integration by parts, we get
\[\sum_{i,j=1}^2\int_{\Omega_L}\p_iu^j\p_ju^i
=-\sum_{i,j=1}^2\int\p_i\p_ju^ju^i+\sum_{i,j=1}^2\int_{\p\Omega_L}\p_iu^ju^in^jdS.\]
Then, it follows from the incompressibility and boundary condition $\mathbf{u}\cdot\mathbf{n}=0$ that the right-hand side of the above equality is 0. Consequence, (\ref{c.9}) follows from \eqref{c.10} together with \eqref{c.3}.
\endProof

In addition, we have
\begin{cor}\label{cc.11}
{\bf(Korn's inequality on $\Omega$)}~~There exists a constant $C>0$, such that for any $\mathbf{u}\in\mathbf{V}$, there holds
\begin{align}\label{c.11}
\|\mathbb{D}(\mathbf{u})\|_{L^2}\geq C\|\mathbf{u}\|_{H^1}.
\end{align}
\end{cor}

\begin{lem}\label{mc.14}
{\bf($L^\infty(\Omega_L)$ estimate)}~~There exists a constant $C>0$, being independent of $L$, such that for any $\mathbf{u}\in \mathcal{W}_L$, there holds
\begin{align}\label{c.16}
\|\mathbf{u}\|^2_{L^\infty(\Omega_L)}\leq C\|\mathbf{u}\|_{L^2(\Omega_L)}\|\mathbf{u}\|_{H^2(\Omega_L)}.
\end{align}
\end{lem}
\Proof
First prove (\ref{c.16}) for $u^1$. Similar to the proof of Lemma \ref{mc.6}, we denote $v^1(x,y)=\zeta(2y)\tilde{u}^1(x,y)$, where
\begin{equation*}
\tilde{u}^1(x,y):=
\begin{cases}
e^{\frac{k_0}{\mu}y}u^1(x,y),&\ y\in [0,1]\\
e^{\frac{-k_0}{\mu}y}u^1(x,-y),&\ y\in [-1,0].
\end{cases}
\end{equation*}
Then, $v^1(x,y)$ satisfies
\[v^1(\pm L,y)=v^1(x,\pm 1)=0,\]
and hence we have
\begin{align*}
|v^1(x,y)|^2&=2\int_{-1}^yv^1_y(x,\theta)v^1(x,\theta)d\theta=2\int_{-L}^x\int_{-1}^y[v^1_{xy}v^1+v^1_yv_x^1](s,\theta)d\theta ds\nonumber\\
&\leq 2\int_{-L}^L\int_{-1}^1|v^1_{xy}v^1+v^1_yv_x^1|(x,y)dxdy.
\end{align*}
As $v^1$ is symmetric in the vertical direction and vanishes on the boundary, we use integration by parts together with H\"older inequality and yield
\begin{align*}
\sup_{(x,y)\in\Omega_L}|v^1(x,y)|^2\leq C\|v^1\|_{L^2(\Omega_L)}\|\nabla^2v^1\|_{L^2(\Omega_L)}.
\end{align*}
Further, in virtue of the definition of $v^1$, we infer that
\begin{align}\label{c.17}
\sup_{(x,y)\in[-L,L]\times[0,1/2]}|u^1(x,y)|^2\leq C\|u^1\|_{L^2(\Omega_L)}\|u^1\|_{H^2(\Omega_L)}.
\end{align}
Similarly, we have
\begin{align}\label{c.18}
\sup_{(x,y)\in[-L,L]\times[1/2,1]}|u^1(x,y)|^2\leq C\|u^1\|_{L^2(\Omega_L)}\|u^1\|_{H^2(\Omega_L)}.
\end{align}
This completes the proof of \eqref{c.16} for $u^1$.

For $u^2$, denote $v^2:=\zeta(\frac{x+L}{L})\tilde{u}^2$, in which
\begin{equation*}
\tilde{u}^2(x,y):=
\begin{cases}
u^2(x,y),&\ x\in [-L,L]\\
u^2(-x-2L,y),&\ x\in [-3L,-L].
\end{cases}
\end{equation*}
Obviously, $v^2(x,y)$ satisfies
\[v^2(-3L,y)=v^2(L,y)=v^2(x,0)=v^2(x,1)=0,\]
and hence, we have
\begin{align*}
|v^2(x,y)|^2\leq 2\int_{-3L}^L\int_0^1|v^2_{xy}v^2+v^2_yv_x^2|(x,y)dxdy.
\end{align*}
Since $v^2$ is symmetric in the horizontal direction and vanishes on the boundary, similarly, we obtain
\begin{align}\label{c.19}
\sup_{(x,y)\in\Omega_L}|v^2(x,y)|^2\leq C\|v^2\|_{L^2(\Omega_L)}\|\nabla^2v^2\|_{L^2(\Omega_L)}.
\end{align}
Further applying the definition of $v^2$ leads to
\begin{align}\label{c.20}
\sup_{(x,y)\in[-L,0]\times[0,1]}|u^1(x,y)|^2\leq C\|u^2\|_{L^2(\Omega_L)}\|u^2\|_{H^2(\Omega_L)}.
\end{align}
Similarly, we also have
\begin{align}\label{c.21}
\sup_{(x,y)\in[0,L]\times[0,1]}|u^2(x,y)|^2\leq C\|u^2\|_{L^2(\Omega_L)}\|u^2\|_{H^2(\Omega_L)}.
\end{align}
The proof of this lemma is finished.
\endProof

As constant $C>0$ in inequality (\ref{c.16}) is independent of $L$, we take $L\rightarrow\infty$ and yield

\begin{cor}\label{cc.15}
{\bf ($L^\infty(\Omega)$ estimate)}~~There exists a constant $C>0$, such that for any $\mathbf{u}\in \mathcal{W}$, it holds
\begin{align}\label{c.22}
\|\mathbf{u}\|^2_{L^\infty}\leq C\|\mathbf{u}\|_{L^2}\|\mathbf{u}\|_{H^2}.
\end{align}
\end{cor}

\section{Proof of Theorem \ref{thm1}}
In this section, we prove the existence and uniqueness of the weak solution to initial boundary value problem (\ref{b.1}) in $\Omega$.

In virtue of Lemma \ref{mc.3}, the function space $V_L$ possess a basis $\{\mathbf{v}_j\}_{j=1}^\infty\subset H^3(\Omega_L)$, which is also a orthonormal basis of $H_L$. For any fixed $m\in\mathbb{N}_+$, we seek approximate solutions in the form $\mathbf{u}_m(t)=\sum_{j=1}^mg^m_j(t)\mathbf{v}_j$, satisfying
\begin{align}\label{c.23}
&\frac{d}{dt}\int_{\Omega_L}\mathbf{u}_m(t)\cdot\mathbf{v}_k+2\mu\int_{\Omega_L}\mathbb{D}(\mathbf{u}_m(t)):\mathbb{D}(\mathbf{v}_k)
+\int_{\Omega_L}\mathbf{u}_m(t)\cdot\nabla\mathbf{u}_m(t)\cdot\mathbf{v}_k\nonumber\\
&=\int_{\p\Omega_L}k(x,y)\mathbf{u}_m(t)\cdot\mathbf{v}_kdS+\int_{\Omega_L}\mathbf{f}\cdot\mathbf{v}_k
\end{align}
for any $k=1,2,\cdots,m$ and the initial data
\begin{align}\label{c.24}
\mathbf{u}_m(0)=\sum_{j=1}^m(\mathbf{u}_0,\mathbf{v}_j)_{V_L}\mathbf{v}_j.
\end{align}
Combining (\ref{c.23}) and (\ref{c.24}) gives a Cauchy problem of ODEs for $(g^m_1(t),g^m_2(t),\cdots,g^m_m(t))$, in which the nonlinear terms is the zeroth-order ones. By the classical theory of the first order ODEs, it possesses a unique solution $(g^m_1(t),g^m_2(t),\cdots,g^m_m(t))\in C^1[0,T_m)$, with $T_m$ being the maximum life time. Hence, there exists a unique solution $\mathbf{u}_m(t)\in C^1([0,T_m);\mathcal{W}_L)$ to problem (\ref{c.23})-(\ref{c.24}). In order to take $m\rightarrow\infty$ and extend $T_m$ to $T$, we need some uniform energy estimates.

{\it Step 1.}~~Multiplying $g^m_k(t)$ to both sides of (\ref{c.23}) and adding them up from $k=1$ to $k=m$, integrating the results by parts yield
\begin{align}\label{c.25}
\frac{1}{2}\frac{d}{dt}\int_{\Omega_L}|\mathbf{u}_m(t)|^2+2\mu\int_{\Omega_L}|\mathbb{D}(\mathbf{u}_m(t))|^2
&=\int_{\p\Omega_L}k(x,y)|\mathbf{u}_m(t)|^2+\int_{\Omega_L}\mathbf{f}\cdot\mathbf{u}_m(t).
\end{align}
Note that
\begin{align}
\int_{\p\Omega_L}k(x,y)|\mathbf{u}_m(t)|^2&=\int^L_{-L}(k_1|u^1(x,1,t)|^2+k_0|u^1(x,0,t)|^2)dx\nonumber\\
&=\int_{\Omega_L}[((k_1+k_0)y-k_0)|u^1(x,y,t)|^2]_ydxdy\nonumber\\
&\leq C\int_{\Omega_L}|\mathbf{u}_m(t)|^2+\e\int_{\Omega_L}|\nabla\mathbf{u}_m(t)|^2,\label{c.26}\\
\int_{\Omega_L}\mathbf{f}\cdot\mathbf{u}_m(t)&\leq \int_{\Omega_L}|\mathbf{u}_m(t)|^2+\int_{\Omega_L}|\mathbf{f}|^2.\label{c.27}
\end{align}
Substituting (\ref{c.26}) and (\ref{c.27}) into (\ref{c.25}), together with using Korn's inequality (\ref{c.9}) and taking $\e$ small sufficiently, one has
\begin{align*}
\frac{d}{dt}\|\mathbf{u}_m(t)\|^2_{L^2(\Omega_L)}+\|\mathbf{u}_m(t)\|^2_{H^1(\Omega_L)}\leq C\|\mathbf{u}_m(t)\|^2_{L^2(\Omega_L)}+C\|\mathbf{f}\|^2_{L^2(\Omega_L)}.
\end{align*}
Then, applying Gronwall's inequality gives
\begin{align}\label{c.28}
\sup_{0\leq t\leq T}\|\mathbf{u}_m(t)\|^2_{L^2(\Omega_L)}+\int_0^T\|\mathbf{u}_m(t)\|^2_{H^1(\Omega_L)}\leq C\left(\|\mathbf{u}_0\|^2_{L^2(\Omega_L)}+\int_0^T\|\mathbf{f}(t)\|^2_{L^2(\Omega_L)}\right).
\end{align}

{\it Step 2.}~~Multiplying both sides of (\ref{c.23}) by $\frac{d}{dt}g^m_k(t)$ and adding them up respect to $k$ from $1$ to $m$, using integration by parts formula, one obtains
\begin{align}\label{c.29}
\mu\frac{d}{dt}\int_{\Omega_L}|\mathbb{D}(\mathbf{u}_m)(t)|^2+\int_{\Omega_L}|\p_t\mathbf{u}_m(t)|^2
\leq&\int_{\Omega_L}|\mathbf{u}_m(t)|^2|\nabla\p_t\mathbf{u}_m(t)|+\int_{\Omega_L}|\mathbf{f}||\p_t\mathbf{u}_m(t)|\nonumber\\
&+\int_{\p\Omega_L}k(x,y)\mathbf{u}_m(t)\cdot\p_t\mathbf{u}_m(t).
\end{align}
Similar to {\it Step 1}, it follows from the Poincar\'e inequality (\ref{c.3}) for $\p_t\mathbf{u}$ that
\begin{align}
\int_{\p\Omega_L}k(x,y)\mathbf{u}_m(t)\cdot\p_t\mathbf{u}_m(t)
&\leq\e\|\nabla\p_t\mathbf{u}_m(t)\|^2_{L^2(\Omega_L)}+\|\mathbf{u}_m(t)\|^2_{H^1(\Omega_L)},\label{c.30}\\
\int_{\Omega_L}|\mathbf{f}(t)||\p_t\mathbf{u}_m(t)|
&\leq \e\|\nabla\p_t\mathbf{u}_m(t)\|^2_{L^2(\Omega_L)}+C\|\mathbf{f}(t)\|^2_{L^2(\Omega_L)}.\label{c.31}
\end{align}
Besides, using H\"older inequality and Lemma \ref{mc.6} deduces
\begin{align}\label{c.32}
\int_{\Omega_L}|\mathbf{u}_m(t)|^2|\nabla\p_t\mathbf{u}_m(t)|
\leq\e\|\nabla\p_t\mathbf{u}_m(t)\|^2_{L^2(\Omega_L)}+\|\mathbf{u}_m(t)\|^2_{L^2(\Omega_L)}\|\mathbf{u}_m(t)\|^2_{H^1(\Omega_L)}.
\end{align}
Substituting (\ref{c.30})-(\ref{c.32}) into (\ref{c.29}) gives
\begin{align}\label{c.33}
\mu\frac{d}{dt}\int_{\Omega_L}|\mathbb{D}(\mathbf{u}_m)(t)|^2+\int_{\Omega_L}|\p_t\mathbf{u}_m(t)|^2
\leq&\e\|\nabla\p_t\mathbf{u}_m(t)\|^2_{L^2(\Omega_L)}+C\|\mathbf{f}(t)\|^2_{L^2(\Omega_L)}\nonumber\\
&+C\left(\|\mathbf{u}_m(t)\|^2_{L^2(\Omega_L)}+1\right)\|\mathbf{u}_m(t)\|^2_{H^1(\Omega_L)}.
\end{align}

Now, applying operator $\frac{d}{dt}$ to (\ref{c.23}) and repeat {\it Step 2}, one gets
\begin{align}\label{c.34}
&\frac{1}{2}\frac{d}{dt}\int_{\Omega_L}|\p_t\mathbf{u}_m(t)|^2+\int_{\Omega_L}|\mathbb{D}(\p_t\mathbf{u}_m)(t)|^2\nonumber\\
\leq&\int_{\Omega_L}|\p_t\mathbf{u}_m(t)|^2|\nabla\mathbf{u}_m(t)|+\int_{\Omega_L}|\p_t\mathbf{f}||\p_t\mathbf{u}_m(t)|
+\int_{\p\Omega_L}k(x,y)|\p_t\mathbf{u}_m(t)|^2\nonumber\\
\leq&\e\|\nabla\p_t\mathbf{u}_m(t)\|^2_{L^2(\Omega_L)}+\|\p_t\mathbf{f}(t)\|^2_{L^2(\Omega_L)}
+\left(\|\mathbf{u}_m(t)\|^2_{H^1(\Omega_L)}+1\right)\|\p_t\mathbf{u}_m(t)\|^2_{L^2(\Omega_L)}.
\end{align}

Adding up (\ref{c.33}) and (\ref{c.34}) with $\e$ small sufficiently, integrating the result over $(0,t)$ and using Korn's inequality (\ref{c.9}), we reach
\begin{align}\label{c.35}
&\|\p_t\mathbf{u}_m(t)\|^2_{L^2(\Omega_L)}+\|\mathbf{u}_m(t)\|^2_{H^1(\Omega_L)}
+\int_0^t\|\p_t\mathbf{u}_m(s)\|^2_{H^1(\Omega_L)}ds\nonumber\\
\leq& \int_0^t\left(\|\mathbf{u}_m(s)\|^2_{H^1(\Omega_L)}+1\right)
\left(\|\p_t\mathbf{u}_m(s)\|^2_{L^2(\Omega_L)}+\|\mathbf{u}_m(s)\|^2_{H^1(\Omega_L)}\right)ds\nonumber\\
&+\|\p_t\mathbf{u}_m(0)\|^2_{L^2(\Omega_L)}+\|\mathbf{u}(0)\|^2_{H^1(\Omega_L)}
+\int_0^t\|(\mathbf{f},\p_t\mathbf{f})(s)\|_{L^2(\Omega_L)}^2ds.
\end{align}
In the next step, we should estimate $\|\p_t\mathbf{u}_m(0)\|^2_{L^2(\Omega_L)}$.

{\it Step 3.} In fact, different from (\ref{c.29}), multiplying (\ref{c.23}) by $\frac{d}{dt}g^m_k(t)$, adding them up from $k=1$ to $k=m$, and using integration by parts formula, we also have
\begin{align*}
\int_{\Omega_L}|\p_t\mathbf{u}_m(t)|^2=&-\mu\int_{\Omega_L}\Delta\mathbf{u}_m(t)\cdot\p_t\mathbf{u}_m(t)
+\int_{\Omega_L}\mathbf{f}\cdot\p_t\mathbf{u}_m(t)\nonumber\\
&+\int_{\Omega_L}\mathbf{u}_m(t)\cdot\nabla\mathbf{u}_m(t)\cdot\p_t\mathbf{u}_m(t)\\
\leq& \frac{3}{4}\|\p_t\mathbf{u}_m(t)\|^2_{L^2(\Omega_L)}+\|\mathbf{u}_m(t)\cdot\nabla\mathbf{u}_m(t)\|^2_{L^2(\Omega_L)}\\
&+\|\Delta\mathbf{u}(t)\|^2_{L^2(\Omega_L)}+\|\mathbf{f}\|^2_{L^2(\Omega_L)}.
\end{align*}
Taking $t\rightarrow0$ and using (\ref{c.16}), we get
\[\|\p_t\mathbf{u}_m(0)\|_{L^2(\Omega_L)}\leq C\left(1+\|\mathbf{u}_0\|^2_{H^2(\Omega_L)}+\|\mathbf{f}(0)\|_{L^2(\Omega_L)}\right).\]

Now, substituting it into (\ref{c.35}) and using Gronwall's inequality together with (\ref{c.28}) give
\begin{align}\label{c.36}
&\sup_{0\leq t\leq T}\left(\|\p_t\mathbf{u}_m(t)\|^2_{L^2(\Omega_L)}+\|\mathbf{u}_m(t)\|^2_{H^1(\Omega_L)}\right)
+\int_0^T\|\p_t\mathbf{u}_m(t)\|^2_{H^1(\Omega_L)}dt\nonumber\\
\leq& C\left(T,\|\mathbf{u}_0\|_{H^2(\Omega_L)},\|\mathbf{f}(0)\|_{L^2(\Omega_L)},\|\mathbf{f}\|_{H^1(0,T;L^2(\Omega_L))}\right).
\end{align}
Besides, it holds that
\[\|\mathbf{f}(0)\|_{L^2(\Omega_L)}^2=\int_0^T\p_t\left(\frac{s-T}{T}\|\mathbf{f}(s)\|_{L^2(\Omega_L)}^2\right)ds
\leq C\left(T,\|\mathbf{f}(t)\|^2_{H^1(0,T;L^2(\Omega_L))}\right).\]
Then, (\ref{c.36}) can be further simplified as
\begin{align*}
&\sup_{0\leq t\leq T}\left(\|\p_t\mathbf{u}_m(t)\|^2_{L^2(\Omega_L)}+\|\mathbf{u}_m(t)\|^2_{H^1(\Omega_L)}\right)
+\int_0^T\|\p_t\mathbf{u}_m(t)\|^2_{H^1(\Omega_L)}dt\nonumber\\
\leq& C\left(T,\|\mathbf{u}_0\|_{H^2(\Omega_L)},\|\mathbf{f}\|_{H^1(0,T;L^2(\Omega_L))}\right).
\end{align*}
{\it Step 4.} This means that the sequence $\{\mathbf{u}_m(t)\}_{m=1}^\infty$ is bounded in $L^\infty(0,T;\mathbf{V}_L)\cap L^2(0,T;\mathbf{V}_L)$, and hence the maximum life time $T_m$ can be extended to $T$. In addition, it also tells that the sequence $ \{\p_t\mathbf{u}_m(t)\}_{m=1}^\infty$ is bounded in $L^\infty(0,T;\mathbf{H}_L)\cap L^2(0,T;\mathbf{V}_L)$. Consequently, by Aubin-Lions Lemma, there exists a function
$\mathbf{u}^L\in C([0,T];\mathbf{H}_L)\cap L^2(0,T;\mathbf{V}_L)$ and a subsequence
 $ \{\mathbf{u}_{m'}(t)\}^\infty_{m'=1}\subset\{\mathbf{u}_m(t)\}^\infty_{m=1}$, such that as $m'\rightarrow\infty$,
\begin{align*}
&\mathbf{u}_{m'}\rightarrow\mathbf{u}^L\ \textrm{ Stongly~in~~}\ C([0,T];\mathbf{H}_L),\\
&\mathbf{u}_{m'}\rightarrow\mathbf{u}^L\ \textrm{weakly*~in~~}\ L^\infty(0,T;\mathbf{V}_L),\\
&\p_t\mathbf{u}_{m'}\rightarrow\p_t\mathbf{u}^L\ \textrm{weakly*~in~~}\ L^\infty(0,T;\mathbf{H}_L).
\end{align*}

Since for any $\mathbf{v}\in\mathbf{V}_L$, the subsequence $\{\mathbf{u}_{m'}(t)\}^\infty_{m'=1}$ satisfies
\begin{align*}
&\frac{d}{dt}\int_{\Omega_L}\mathbf{u}_{m'}(t)\cdot\mathbf{v}+2\mu\int_{\Omega_L}\mathbb{D}(\mathbf{u}_{m'}(t)):\mathbb{D}(\mathbf{v})
+\int_{\Omega_L}\mathbf{u}_{m'}(t)\cdot\nabla\mathbf{u}_{m'}(t)\cdot\mathbf{v}\\
&=\int_{\p\Omega_L}k(x,y)\mathbf{u}_{m'}(t)\cdot\mathbf{v}dS+\int_{\Omega_L}\mathbf{f}\cdot\mathbf{v},
\end{align*}
where the first term in the right-hand side can be rewritten as
\begin{align}\label{c.37}
\int_{\Omega_L}[((k_1+k_0)y-k_0)u^1_{m'}v^1]_y
=&(k_1+k_0)\int_{\Omega_L}u^1_{m'}v^1+\int_{\Omega_L}[(k_1+k_0)y-k_0]\p_yu^1_{m'}v^1\nonumber\\
&+\int_{\Omega_L}[(k_1+k_0)y-k_0]u^1_{m'}\p_yv^1.
\end{align}
Thus, taking $m'\rightarrow\infty$ deduces that $\mathbf{u}^L$ satisfies
\begin{align}\label{c.38}
&\frac{d}{dt}\int_{\Omega_L}\mathbf{u}^L(t)\cdot\mathbf{v}+2\mu\int_{\Omega_L}\mathbb{D}(\mathbf{u}^L(t)):\mathbb{D}(\mathbf{v})
+\int_{\Omega_L}\mathbf{u}(t)\cdot\nabla\mathbf{u}^L(t)\cdot\mathbf{v}\nonumber\\
&=\int_{\p\Omega_L}k(x,y)\mathbf{u}^L(t)\cdot\mathbf{v}dS+\int_{\Omega_L}\mathbf{f}\cdot\mathbf{v}, \ \textrm{for~any~} \mathbf{v}\in\mathbf{V}_L,
\end{align}
Correspondingly, by weak lower continuity\cite{Temam84}, there holds
\begin{align}\label{c.39}
&\sup_{0\leq t\leq T}\left(\|\mathbf{u}^L(t)\|^2_{H^1(\Omega_L)}+\|\p_t\mathbf{u}^L(t)\|^2_{H^1(\Omega_L)}\right)
+\int_0^T\left(\|\mathbf{u}^L(t)\|^2_{H^1(\Omega_L)}+\|\p_t\mathbf{u}^L(t)\|^2_{H^1(\Omega_L)}\right)\nonumber\\
&\leq C\left(T,\|\mathbf{u}_0\|_{H^2(\Omega_L)},\|\mathbf{f}\|_{H^1(0,T;L^2(\Omega_L))}\right).
\end{align}
Consequently, one has
\[\mathbf{u}^L\in L^\infty(0,T;\mathbf{V}_L)\cap L^2(0,T;\mathbf{V}_L),\ \p_t\mathbf{u}^L\in L^\infty(0,T;\mathbf{H}_L)\cap L^2(0,T;\mathbf{V}_L).\]
Moreover, the estimates for the weak solutions are independent of $L$, and hence the existence of the weak solution and estimate (\ref{c.2}) follows so long as taking $L\rightarrow\infty$ .

{\it Step 5.} In this step, we prove the uniqueness of the weak solution. Assume that there are two weak solutions $\mathbf{u}_1(t),\mathbf{u}_2(t)\in L^\infty(0,T;\mathbf{V})\cap L^2(0,T;\mathbf{V})$  to problem (\ref{b.1}), satisfying weak formulation (\ref{c.38}), estimate (\ref{c.39}) and the initial data
$\mathbf{u}_1(0) =\mathbf{u}_2(0)=\mathbf{u}_0.$ Then, the difference $\bar{\mathbf{u}}(t):=\mathbf{u}_1(t)-\mathbf{u}_2(t)$ satisfies weak formula
\begin{align}\label{c.40}
&\frac{d}{dt}\int_\Omega\bar{\mathbf{u}}(t)\cdot\mathbf{v}+2\mu\int_\Omega\mathbb{D}(\bar{\mathbf{u}}(t)):\mathbb{D}(\mathbf{v})
+\int_\Omega\mathbf{u}_1(t)\cdot\nabla\bar{\mathbf{u}}(t)\cdot\mathbf{v}
+\int_\Omega\bar{\mathbf{u}}\cdot\nabla\mathbf{u}_2(t)\cdot\mathbf{v}\nonumber\\
&=\int_{\p\Omega}k(x,y)\bar{\mathbf{u}}(t)\cdot\mathbf{v}dS.
\end{align}

Specially take $\mathbf{v}=\bar{\mathbf{u}}$. Then, integrating by parts, we reach
\begin{align}\label{c.41}
\frac{d}{dt}\|\bar{\mathbf{u}}(t)\|_{L^2}^2+\|\bar{\mathbf{u}}(t)\|^2_{H^1}
\leq C\int_\Omega|\bar{\mathbf{u}}|^2|\nabla\mathbf{u}_2(t)|+C\int_{\p\Omega}k(x,y)|\bar{u}^1(t)|^2dS.
\end{align}
Using the skill in (\ref{c.37}), we find that
\begin{align}\label{c.42}
\int_{\p\Omega}k(x,y)|\bar{u}^1(t)|^2dS&\leq C\|\bar{\mathbf{u}}(t)\|^2_{L^2}+C\|\bar{\mathbf{u}}(t)\|_{L^2}\|\bar{\mathbf{u}}(t)\|_{H^1}\nonumber\\
&\leq \e\|\bar{\mathbf{u}}(t)\|_{H^1}+ C\|\bar{\mathbf{u}}(t)\|_{L^2}.
\end{align}
In addition, it follows from (\ref{c.6}) that
\begin{align}\label{c.43}
\int_\Omega|\bar{\mathbf{u}}|^2|\nabla\mathbf{u}_2(t)|&\leq \|\bar{\mathbf{u}}(t)\|^2_{L^4}\|\nabla\mathbf{u}_2(t)\|_{L^2}\nonumber\\
&\leq \e\|\bar{\mathbf{u}}(t)\|_{H^1}+ C\|\bar{\mathbf{u}}(t)\|^2_{L^2}\|\nabla\mathbf{u}_2(t)\|^2_{L^2}
\end{align}

Then, substituting (\ref{c.42}) and (\ref{c.43}) into (\ref{c.41}) with $\e$ being small sufficiently implies
\begin{align}\label{c.44}
\frac{d}{dt}\|\bar{\mathbf{u}}(t)\|_{L^2}^2\leq C\left(1+\|\mathbf{u}_2(t)\|^2_{H^1}\right)\|\bar{\mathbf{u}}(t)\|^2_{L^2}
\end{align}
Finally, applying Gronwall's inequality to (\ref{c.44}) with estimate(\ref{c.39}) and the fact that $\bar{\mathbf{u}}(0)\equiv 0$, we yield $\bar{\mathbf{u}}\equiv 0$, which completes the proof of uniqueness. \endProof

\section{Proof of Theorem \ref{thm2}}
In order to obtain higher order energy estimates and thus imply that the weak solution is in fact a strong solution, even smooth solution, we need the following Stokes estimate.

\begin{prop}\label{pc.16}
Assume that $\mathbf{u}\in H^1$ is the weak solution to the following initial boundary value problem
\begin{align}\label{c.45}
\begin{cases}
-\mu\Delta\mathbf{u}+\nabla p=\mathbf{F},&\ \textrm{in~~}\Omega,\\
\nabla\cdot\mathbf{u}=0,&\ \textrm{in~~}\Omega,\\
\mathbf{u}\cdot\mathbf{n}=0,&\ \textrm{on~~}\p\Omega,\\
2\mu\mathbb{D}(\mathbf{u})\mathbf{n}\cdot\tau=k(x,y)\mathbf{u}\cdot\tau,&\ \textrm{on~~}\p\Omega,
\end{cases}
\end{align}
where $\mathbf{F}\in L^2$ and $k(x,1)=k_1,k(x,0)=k_0$ are constants. Then $\mathbf{u}\in H^2$ and satisfies
\begin{align}\label{c.46}
\|\mathbf{u}\|_{H^2}+\|\nabla p\|_{L^2}\leq C\left(\|\mathbf{F}\|_{L^2}+\|\mathbf{u}\|_{L^2}\right),
\end{align}
where $C>0$ depends only on $\mu,k_0,k_1. $
\end{prop}

\Proof The proof of this proposition consists of 4 steps.

 {\it Step 1.}~~For any positive constant $ \beta$ large enough and functions $\mathbf{u},\mathbf{F},k(x,y)$ given in (\ref{c.45}), we consider the auxiliary problem:
\begin{align}\label{c.47}
\begin{cases}
-\mu\Delta w+\beta w=\mathrm{curl}\mathbf{F}+\beta\mathrm{curl}\mathbf{u}:=\mathrm{curl}\mathbf{\Phi},&\ \textrm{in~~}\Omega,\\
w=k(x,y)\mathbf{u}\cdot\tau:=g,&\ \textrm{on~~}\p\Omega.
\end{cases}
\end{align}

Since $\mathrm{curl}\mathbf{\Phi}\in H^{-1}$, we define the bilinear form
\begin{align}\label{c.48}
\mathscr{B}[w,\tilde{w}]=\mu\int_{\Omega}\nabla w:\nabla\tilde{w}+\beta\int_{\Omega}w\tilde{w},
\end{align}
 for $w,\tilde{w}\in H^1_g:=\{w\in H^1|w=g\ \textrm{in~~}\p\Omega\}$. As the inhomogeneous Dirichlet boundary value problem problem \eqref{c.47} can be rewritten as a homogeneous one via homogenization method, without lose of any generality, we assume that $g=0$. It is easy to check that the bilinear form $\mathscr{B}$ is continuous and coercive on $H^1_0$, and hence it follows from the Lax-Milgram theorem that there exists an unique $w\in H^1_g$ being the weak solution to system (\ref{c.47}), i.e.
\begin{align}\label{c.49}
\mu\int_{\Omega}\nabla w\cdot\nabla\tilde{w}+\beta\int_\Omega w\tilde{w}=-\int_\Omega\mathbf{\Phi}\cdot\overrightarrow{\mathrm{curl}}\tilde{w},
\end{align}
holds for any $\tilde{w}\in H^1_0$. Here $\overrightarrow{\mathrm{curl}}:=(-\p_y,\p_x)$.

Take $\tilde{w}=w-[(k_0+k_1)y-k_0]u^1$. Then $\tilde{w}\in H^1_0$. Substituting it into (\ref{c.49}) and using Cauchy problem gives
\begin{align*}
\mu\int_\Omega|\nabla  w|^2+\beta\int_\Omega|w|^2\leq \frac{\mu}{2}\int_\Omega|\nabla w|^2+\frac{\beta}{2}\int_\Omega|w|^2+C\int_\Omega\left(|\mathbf{\Phi}|^2+|\mathbf{u}|^2+|\nabla \mathbf{u}|^2\right).
\end{align*}
This indicates that
\begin{align}\label{c.50}
\|w\|_{H^1}\leq C\left(\|\mathbf{\Phi}\|_{L^2}+\|\mathbf{u}\|_{H^1}\right).
\end{align}
{\it Step 2.}~~For $w$ constructed in {\it Step 1}, consider the following boundary value problem
\begin{align}\label{c.51}
\begin{cases}
-\Delta\Psi=-w,&\ \textrm{in~~}\Omega,\\
\Psi=0,&\ \textrm{on~~}\p\Omega.
\end{cases}
\end{align}
By the classical elliptic equation theory, problem (\ref{c.51}) possesses a unique solution $\Psi\in H^3$. In what follows, we deduce the $H^3$ estimate for $\Psi$.

Multiplying (\ref{c.51})$_1$ by $\Psi$, integrating by parts over $\Omega$ and using Poincar\'e inequality, we get
\[ \|\Psi\|_{H^1}\leq C\|w\|_{L^2}.\]
Applying $\p_x$ to (\ref{c.51})$_1$, similarly, we deduce
\[\|\Psi_x\|_{H^1}\leq C\|w\|_{L^2}.\]
Moreover, it follows from (\ref{c.50})$_1$ that $ \Psi_{yy}=w-\Psi_{xx}$. Then we also have \[\|\Psi_{yy}\|_{L^2}\leq \|w\|_{L^2}+\|\Psi_{xx}\|_{L^2}\leq C\|w\|_{L^2}.\]
In conclusion, we yield $\|\Psi\|_{H^2}\leq C\|w\|_{L^2}$.

Note again that $\Psi_x$ satisfies
\begin{align}\label{c.52}
\begin{cases}
-\Delta\Psi_x=-w_x,&\ \textrm{in~~}\Omega,\\
\Psi_x=0,&\ \textrm{on~~}\p\Omega.
\end{cases}
\end{align}
Then, by the analysis above, one has $\|\Psi_x\|_{H^2}\leq C\|w_x\|_{L^2}$. To obtain estimates for $\Psi_{yyy}$, we apply $\p_y$  to (\ref{c.51})$_1$ and yield $\Psi_{yyy}=\p_y w-\p_y\Psi_{xx}$, which leads to
\[\|\Psi_{yyy}\|_{L^2}\leq C\|w\|_{H^1}+\|\Psi_x\|_{H^2}\leq C\|w\|_{H^1}.\]
Thus, there holds $\|\Psi\|_{H^3}\leq C\|w\|_{H^1}$.

{\it Step 3.}~~Go back to problem (\ref{c.45}) and take $\mathbf{v}=\overrightarrow{\mathrm{curl}}\Psi$. Then $\mathbf{v}\in H^2$ satisfies
\begin{align}\label{c.53}
\|\mathbf{v}\|_{H^2}\leq C\|w\|_{H^1}
\end{align}
and the relationship $w=\mathrm{curl}\mathbf{v}$. Furthermore, substituting this relationship into (\ref{c.47}), we reach
\begin{align}\label{c.54}
\begin{cases}
-\mu\Delta\mathrm{curl}\mathbf{v}+\beta\mathrm{curl}\mathbf{v}=\mathrm{curl}(\mathbf{F}+\beta\mathbf{u}),&\ \textrm{in~~}\Omega,\\
\mathrm{curl}\mathbf{v}=k(x,y)\mathbf{u}\cdot\tau,&\ \textrm{on~~}\p\Omega,\\
\end{cases}
\end{align}
In virtue of Hodge decomposition, equation (\ref{c.54})$_1$ is equivalent to
\[-\mu\Delta\mathbf{v}+\beta\mathbf{v}+\nabla q=\mathbf{F}+\beta\mathbf{u},\ \textrm{in~~}\Omega.\]
In the meantime, by Lemma \ref{mc.2}, boundary condition (\ref{c.54})$_2$ is equivalent to \[2\mathbb{D}(\mathbf{v})\mathbf{n}\cdot\tau=k(x,y)\mathbf{u}\cdot\tau,\ \textrm{on~~}\p\Omega.\] Besides, it follows from the definition of $\mathbf{v}$ that $\mathrm{div}\mathbf{v}=0.$ In addition, since $ \Psi|_{\p\Omega}=0$, we have $\nabla\Psi\cdot\tau|_{\p\Omega}=0$, which is equivalent to
\[\overrightarrow{\mathrm{curl}}\Psi\cdot\mathbf{n}|_{\p\Omega}=0,
~~\textrm{i.e.}~~\mathbf{v}\cdot\mathbf{n}|_{\p\Omega}=0.\]

In conclusion, $\mathbf{v}$ is a solution of problem
\begin{align}\label{c.55}
 \begin{cases}
 -\mu\Delta\mathbf{v}+\beta\mathbf{v}=\mathbf{F}+\beta\mathbf{u},&\ \textrm{in~~}\Omega,\\
 \mathrm{div}\mathbf{v}=0,&\ \textrm{in~~}\Omega,\\
 \mathbf{v}\cdot\mathbf{n}=0,&\ \textrm{on~~}\p\Omega\\
 2\mathbb{D}(\mathbf{v})\mathbf{n}\cdot\tau=k(x,y)\mathbf{u}\cdot\tau,&\ \textrm{on~~}\p\Omega.\\
 \end{cases}
\end{align}
{\it Step 4.}~~(\ref{c.45}) indicates that $\mathbf{u}$ is also a weak solution to (\ref{c.55}). Thus, $\mathbf{u}=\mathbf{v}$ is true if the solution of problem (\ref{c.55}) is unique. In fact, one can see that $\mathbf{u}-\mathbf{v}$ satisfies
\begin{align*}
\begin{cases}
 -\mu\Delta(\mathbf{u}-\mathbf{v})+\beta(\mathbf{u}-\mathbf{v})+\nabla(p-q)=0,&\ \textrm{in~~}\Omega,\\
 \mathrm{div}(\mathbf{u}-\mathbf{v})=0,&\ \textrm{in~~}\Omega,\\
 (\mathbf{u}-\mathbf{v})\cdot\mathbf{n}=0,&\ \textrm{on~~}\p\Omega\\
 2\mathbb{D}(\mathbf{u}-\mathbf{v})\mathbf{n}\cdot\tau=0,&\ \textrm{on~~}\p\Omega.\\
 \end{cases}
\end{align*}
It can be deduced from the standard energy method that $\mathbf{u}-\mathbf{v}\equiv 0$, i.e. $\mathbf{u}=\mathbf{v}\in H^2$. Consequently, in virtue of (\ref{c.53}) and (\ref{c.50}), we infer that
\begin{align}\label{c.56}
\|\mathbf{u}\|_{H^2}\leq C\left(\|\mathbf{F}\|_{L^2}+\|\mathbf{u}\|_{H^1}\right).
\end{align}
In particular, substituting (\ref{c.8}) into (\ref{c.56}) and using Cauchy inequality, we get
\begin{align}\label{c.57}
\|\mathbf{u}\|_{H^2}\leq C\left(\|\mathbf{F}\|_{L^2}+\|\mathbf{u}\|_{L^2}\right).
\end{align}
The final work is to deduce estimate for $\nabla p$, which can be directly implied by (\ref{c.45})$_1$ and (\ref{c.57}).

The proof of this proposition is completed.
\endProof

With this Stokes estimate in hand, we are able to state and prove the regularity of the solution.

\ProofTheorem~~~~
By Theorem \ref{thm1}, the initial boundary problem (\ref{b.1}) has a unique weak solution $\mathbf{u}(t)\in L^\infty(0,T;\mathbf{V})$ satisfies estimate (\ref{c.2}). Thus, we still need to prove estimate (\ref{b.2}). In fact, the initial boundary value problem (\ref{b.1}) can be rewritten as
\begin{align}\label{c.59}
\begin{cases}
-\mu\Delta\mathbf{u}+\nabla p=-\p_t\mathbf{u}-\mathbf{u}\cdot\nabla\mathbf{u}+\mathbf{f},&\ \textrm{in~~}\Omega,\\
\nabla\cdot\mathbf{u}=0,&\ \textrm{in~~}\Omega,\\
\mathbf{u}\cdot\mathbf{n}=0,&\ \textrm{on~~}\p\Omega,\\
2\mu\mathbb{D}(\mathbf{u})\mathbf{n}\cdot\tau=k(x,y)\mathbf{u}\cdot\tau,&\ \textrm{on~~}\p\Omega.
\end{cases}
\end{align}
Then, it follows from proposition \ref{pc.16} that $\mathbf{u}(t)\in H^2$ for a.e.$t\in[0,T]$, and that
\begin{align}\label{c.60}
\|\mathbf{u}(t)\|_{H^2}+\|\nabla p(t)\|_{L^2}
\lesssim \|\mathbf{f}(t)\|_{L^2}
+\|\p_t\mathbf{u}\|_{L^2}+\|\mathbf{u}(t)\cdot\nabla\mathbf{u}(t)\|_{L^2}+\|\mathbf{u}\|_{L^2}.
\end{align}
In addition, applying (\ref{c.22}) gives
\begin{align}\label{c.61}
\|\mathbf{u}(t)\cdot\nabla\mathbf{u}(t)\|_{L^2}&\leq \|\mathbf{u}(t)\|_{L^\infty}\|\nabla\mathbf{u}(t)\|_{L^2}\nonumber\\
&\leq \|\mathbf{u}(t)\|^{1/2}_{L^2}\|\mathbf{u}(t)\|^{1/2}_{H^2}\|\nabla\mathbf{u}(t)\|_{L^2}\nonumber\\
&\leq \e\|\mathbf{u}(t)\|_{H^2}+\|\mathbf{u}(t)\|^3_{H^1}.
\end{align}
The proof of this theorem is completed as long as substituting (\ref{c.61}) into (\ref{c.60}) with $\e$ small enough and using estimate (\ref{c.2}).
\endProof

\acknowledgment
The authors would like to thank Prof. Yan Guo from Brown University for the invaluable comments on the present paper and Li also want to express sincere appreciation to Prof. Guo for his kindly assist and academic direction during the period of Li's visitation to Brown University. Ding's research is supported by the National Natural Science Foundation of China (No.11371152, No.11571117, No.11871005 and No.11771155) and Guangdong Provincial Natural Science Foundation (No.2017A030313003).

\end{document}